\input amstex
\documentstyle{amsppt}
\document

\def\irreducible{1.1}
\def\weakirrd{1.2}
\def\symmtensor{1.3}
\def\Zpair{1.4}
\def\Hpairdef{1.5}
\def\hypergraph{2.1}
\def\multihypergraph{2.2}
\def\uniform{2.3}
\def\completedef{2.4}
\def\regular{2.5}
\def\connected{2.6}
\def\partite{2.7}
\def\adten{2.8}
\def\lemmairred{3.1}
\def\Hpair{3.2}
\def\exist{3.3}
\def\redgraph{3.4}
\def\compa{3.5}
\def\compb{3.6}
\def\positive{3.7}
\def\complete{3.8}
\def\compute{3.9}
\def\largestgraph{3.10}
\def\largestcomplete{3.11}
\def\symmetric{3.12}
\def\exampleAa{4.1}
\def\exampleC{4.2}
\def\exampleEe{4.3}
\def\exampleD{4.4}
\def\exampleB{4.5}
\def\exampleA{4.6}
\def\examplecomplete{4.7}
\def\questionA{4.8}
\def\questionB{4.9}
\def\questionC{4.10}
\def\refBerge {1}
\def\refBP {2}
\def\refCS {3}
\def\refCPZ {4}
\def\refCPZPF {5}
\def\refCPZPrim {6}
\def\refCPZZ {7}
\def\refCQZ {8}
\def\refChung {9}
\def\refCocoa {10}
\def\refCooper {11}
\def\refFGH {12}
\def\refHQ {13}
\def\refHu {14}
\def\refKM {15}
\def\refLima {16}
\def\refLimb {17}
\def\refLZI {18}
\def\refNQZ {19}
\def\refQi {20}
\def\refQZ {21}

\hsize=6.5truein
\vsize=9truein

\baselineskip=20pt plus 2pt minus 1.5pt
\NoBlackBoxes

\def\cocoa{{\hbox{\rm C\kern-.13em o\kern-.07em C\kern-.13em o\kern-.15em A}}}

\newcount\qcts\newcount\qcta\qcta=1\newcount\qcteq\newcount\qcthead
\def\nmonth{\ifcase\month\ \or January\or
   February\or March\or April\or May\or June\or July\or August\or
   September\or October\or November\else December\fi}
\def\setref#1{\global\qcteq=0\global\advance\qcts by\number\qcta

\immediate\write3{\string\def\string#1\string{\qctS.\number\qcts\string}}}
\def\seteqn#1{\global\advance\qcteq by\number\qcta\immediate\write3{
    \string\def\string#1\string{\qctS.\number\qcts.\eqnyb\string}}}
\def\ybkey#1{\key{#1}\global\advance\qct by\number\qcta
       \immediate\write3{\string\def\string#1\string{\number\qct\string}}}
\def\eqnyb{\ifnum\qcteq=1 a\else\ifnum\qcteq=2 b\else
   \ifnum\qcteq=3 c\else\ifnum\qcteq=4 d\else
   \ifnum\qcteq=5 e\else\ifnum\qcteq=6 f\else\ifnum\qcteq=7 g\else
   \ifnum\qcteq=8 h\else\ifnum\qcteq=9 i\else
   \ifnum\qcteq=10 j\else\ifnum\qcteq=11 k\else\ifnum\qcteq=12 l\else
   \ifnum\qcteq=13 m\else\ifnum\qcteq=14 n\else\ifnum\qcteq=15 o\else
   \ifnum\qcteq=16 p\else\ifnum\qcteq=17 q\else
   \ifnum\qcteq=18 r\else\ifnum\qcteq=19 s\else\ifnum\qcteq=20 t\else
   \ifnum\qcteq=21u\else\ifnum\qcteq=22 v\else
   \ifnum\qcteq=23 w\else\ifnum\qcteq=24 x\else
   \ifnum\qcteq=25 y\else\ifnum\qcteq=26 z\else *

\fi\fi\fi\fi\fi\fi\fi\fi\fi\fi\fi\fi\fi\fi\fi\fi\fi\fi\fi\fi\fi\fi\fi\fi\fi\fi}
\newcount\qcts\newcount\qcta\qcta=1\newcount\qcteq\newcount\qct

\newcount\qct\newcount\qcta\qct=0\qcta=1
\def\refkey#1{\key{#1}\global\advance\qct
       by\number\qcta
       \immediate\write3{\string\def\string#1
       \string{\number\qct\string}}}

\topmatter
\title On spectral hypergraph theory of the adjacency tensor
\endtitle

\author{Kelly J.~Pearson
\address{Department of Mathematics and Statistics,
         Murray State University, Murray,
         KY 42071-0009.
         \ \ Telephone number: (270) 809-3554 //Department of Mathematics and Statistics,
         Murray State University, Murray,
         KY 42071-0009}
         \endaddress
\email{kpearson\@murraystate.edu}\endemail}
{Tan Zhang
\address{Department of Mathematics and Statistics,
         Murray State University, Murray,
         KY 42071-0009.
         \ \ Telephone number: (270) 809-3554 //Department of Mathematics and Statistics,
         Murray State University, Murray,
         KY 42071-0009}
         \endaddress
\email{tzhang\@murraystate.edu}\endemail}\endauthor


\rightheadtext{} \leftheadtext{} \pageno=1 \abstract  We study both $H$ and $E/Z$-eigenvalues of the adjacency tensor of a uniform multi-hypergraph and give conditions for which the largest positive $H$ or $Z$-eigenvalue corresponds to a strictly positive eigenvector. We also investigate when the $E$-spectrum of the adjacency tensor is symmetric. \endabstract
\endtopmatter

\head{0. Introduction}\endhead
\def\qctS{0}\qcts=0

The spectral theory of graphs emerged as an active research frontier in graph theory since the 1950s. It has since had a major impact on combinatorics, computer science, operations research, biology, and social science. One of the main motivations behind spectral graph theory is to establish connections of the graph's intrinsic structures, such as connectivity, diameters, embeddability, chromatic numbers, with the spectra of various associated matrices, such as the adjacency matrix, the incidence matrix, and the graph Laplacian matrix. More in-depth discussions on spectral graph theory can be found in the seminal book \cite\refChung.

Contrasting the spectral theory of graphs, the spectral theory of hypergraphs is still in its infancy. Some recent development \cite{\refBP, \refCooper, \refHQ, etc.} in this direction has revealed interesting as well as somewhat elusive connections of combinatorial and geometric structures of a given hypergraph,  similar to those of a undirected graph, with the eigen-pairs of special higher-order nonnegative symmetric tensors, such as the adjacency tensor or the Laplacian tensor. However, due to the nonlinear nature of the eigenvalue problems for tensors, different notions of eigenvalues used in the study of a given hypergraph may lead to surprisingly different yet useful conclusions.

Definitions used in this introduction can be found organized in \S 1 Basic Tensor Definitions and \S 2 Basic Hypergraph Definitions.  
The purpose of this paper is to study the $H$ and $E/Z$-eigenvalues of the adjacency tensor of a uniform multi-hypergraph.  

In \cite\refCooper, the authors derive some properties of the eigenvalue of an $m$-graph $\Cal H$ with largest modulus.  In general, there may be many such eigenvalues with the same modulus.  The authors show there is always a real, positive $H$-eigenvalue with maximal modulus.  Moreover, the authors show that any positive $H$-eigenvalue with an $H$-eigenvector whose coordinates are all positive dominates any other positive $H$-eigenvalue with nonnegative $H$-eigenvector, cf. Lemma 3.3, Corollary 3.4,  Lemma 3.5, Corollary 3.6, and Theorem 3.7 of \cite\refCooper. \ After we establish the needed definitions, at the beginning of in \S 3, we summarize their results and give a different proof based on weak irreducibility of tensors, which is first introduced by \cite\refFGH.  The proof given in this paper extends to an $m$-multigraph.

In \cite\refBP, they provided new bounds on the clique number of an undirected graph G, using the $H$-eigenvalue spectral radius of the adjacency tensor of a $k$-clique $(k+1)$-graph of G.

The main results of this paper are given in \S 3 and concern the existence of the largest positive $Z$-eigenvalue and the symmetry of the $E$-spectrum.  For a connected $m$-multigraph on $n$ vertices, we prove the existence of a positive $Z$-eigenvalue with a nonnegative $Z$-eigenvector.  By \cite\refCPZZ, this $Z$-eigenvalue maximizes ${\Cal A}_H x^m$ on the sphere $S^{n-1}$.  Furthermore, if this is a nicely-connected $m$-multigraph, a new concept introduced in this paper, there is a strictly positive eigenvector associated to this eigenvalue.  Bounds for this $Z$-eigenvalue in terms of hypergraph properties are also given.
Examples are given in \S 4 along with some open questions.

\head {1.  Basic Tensor Definitions}\endhead
\def\qctS{1}\qcts=0
We begin this paper with some of the basic definitions from higher order multi-dimensional tensors.  These are standard definitions used in \cite{\refCPZ, \refCPZPF, \refCQZ, \refLima, \refLimb, \refNQZ, \refQi} to name just a few sources from the literature.

Let ${\Bbb R}$ be the real field; we consider an $m$-order
$n$ dimensional tensor ${\Cal A}$ consisting of $n^m$ entries in
${\Bbb R}$:
$${\Cal A}=(a_{i_1\cdots i_m}),\quad  a_{i_1\cdots i_m}\in {\Bbb R}, \quad 1\le
i_1,\ldots, i_m\le n.$$
To an $n$-vector $x=(x_1, \cdots, x_n)$, real or complex, we define
a $n$-vector:
$${\Cal A} x^{m-1}:=\bigg(\sum_{i_2, \ldots, i_m=1}^n a_{i i_2\cdots i_m} x_{i_2}\cdots x_{i_m}\bigg)_{1\le i\le n}.$$
We also define a homogeneous polynomial of degree $m$ in $n$ indeterminants by:
$${\Cal A} x^{m}:=\bigg(\sum_{i_1, \ldots, i_m=1}^n a_{i_1\cdots i_m} x_{i_1}\cdots x_{i_m}\bigg)_{1\le i\le n}.$$
 We denote the set of $m$-order $n$ dimensional tensors via ${\Bbb R}^{[m,n]}$ and the set of $m$-order $n$ dimensional nonnegative tensors via ${\Bbb R}^{[m,n]}_+$ respectively.
\smallbreak
Throughout the remainder of this paper, we shall denote the standard positive cone in $\Bbb R^n$ by ${\Bbb R}_+^n$, namely, ${\Bbb R}_+^n=\{x=(x_1,\cdots x_n)\in {\Bbb R}^n\ |\ x_i\ge 0, 1\le i\le n\}$. We denote the interior of $\Bbb R^n_+$ by $\Bbb R^n_{++}$. 

The following definition has become standard in practice. 
\setref\irreducible
\definition{Definition \irreducible} A tensor $\Cal A=(a_{i_1\cdots i_m})\in \Bbb R^{[m,n]}$  is called reducible, if  there exists a nonempty proper index subset $I\subset \{1, \ldots, n\}$ such that $$a_{i_1\cdots i_m}=0, \quad \forall i_1\in I, \quad \forall i_2, \ldots, i_m\notin I. $$ 
A tensor is called irreducible if it is not reducible.
\enddefinition

The following definition was first introduced by Friedland~et.~al. \cite\refFGH. \ For a given $\Cal A=(a_{i_1\cdots i_m})\in \Bbb R_+^{[m,n]}$, it is
associated to a directed graph $G(\Cal A)=(V, E(\Cal A))$, where
$V=\{1,2,\cdots, n\}$ and a directed edge $(i,j)\in E(\Cal A)$ if there
exists indices $\{i_2,\cdots, i_m\}$ such that $j \in \{i_2,\cdots,
i_m\}$ and $a_{i i_2 \cdots i_m}>0$; in particular, we have
$\sum_{j\in \{i_2,\cdots, i_m \}} a_{ii_2,\cdots i_m }>0$.  A graph is strongly connected  if it contains a directed path from $i$ to $j$ and a directed path from $j$ to $i$ for every pair of vertices $i, j$.
\setref\weakirrd
\definition{Definition \weakirrd} 
A nonnegative tensor $\Cal A\in \Bbb R^{[m,n]}_+$ is called weakly irreducible if the associate directed graph $G(\Cal A)$ is strongly connected.
\enddefinition
\smallbreak
It is equivalent to say (according to  \cite{\refHu}) that the matrix $M(\Cal A)=(m_{ij})$ is irreducible, where
$$ m_{ij}=\Sigma_{j\in \{i_2,\cdots, i_m \}} a_{ii_2,\cdots i_m }.$$

\setref\symmtensor
\definition{Definition \symmtensor}  A tensor $\Cal A \in {\Bbb R}^{[m,n]}$ is symmetric if $a_{i_1 \cdots i_m} = a_{\sigma(i_1 \cdots i_m)}$ for all $\sigma \in \Sigma_m$, the symmetric group on $m$ indices.
\enddefinition

\setref\Zpair
\definition{Definition \Zpair} Let ${\Cal A}$ be a nonzero tensor.   A pair
$(\lambda, x)\in \Bbb{C}\times(\Bbb{C}^n\setminus\{0\})$ is called an $E$-eigenvalue and $E$-eigenvector (or simply $E$-eigenpair) of $\Cal A$  if they satisfy the equations 

$$\eqalign{\Cal A x^{m-1}&=\lambda x \cr
x_1^2+\ldots + x_n^2&=1}$$
We call $(\lambda, x)$ a $Z$-eigenpair if they are both real.
\enddefinition
The $Z$-eigenvalue problem for tensors involves finding nontrivial solutions of inhomogeneous polynomial systems in
several variables.  We define the $Z$-spectrum of $\Cal A$, denoted $\Cal Z(\Cal A)$ to be the set of all $Z$-eigenvalues of $\Cal A$.  It is proven in \cite\refCS, that for a symmetric tensor $\Cal A$, the set of $E$-eigenvalues of $\Cal A$ is nonempty and finite. We therefore define the $Z$-spectral radius of a symmetric tensor $\Cal A$, denoted $\varrho(\Cal A)$, to be $\varrho(\Cal A):=\max\{|\lambda|\ |\ \lambda\in \Cal Z(\Cal A)\}$.  We refer the interested readers to \cite\refCPZZ \ for a more detailed discussion on this topic.

\setref\Hpairdef
\definition{Definition \Hpairdef} Let ${\Cal A}$ be a nonzero tensor.   A pair
$(\lambda, x)\in \Bbb{C}\times(\Bbb{C}^n\setminus\{0\})$ is called an eigenvalue and eigenvector (or simply eigenpair) of $\Cal A$  if they satisfy the equations 
$${\Cal A x^{m-1}=\lambda x^{[m-1]},}$$
where $x^{[m-1]}_i = x_i^{m-1}$.
We call $(\lambda, x)$ a $H$-eigenpair if they are both real.\enddefinition

The $H$-eigenvalue problem for tensors involves finding nontrivial solutions of homogeneous polynomial systems in
several variables. In contrast, the set of eigenvalues of any nonzero tensor is always finite; in which case, the $H$-spectral radius of $\Cal A$, denoted $\rho(\Cal A)$ is defined to be the eigenvalue of $\Cal A$ with maximum modulus.

Previously, both Chang~et.~al. \cite\refCPZPF\ and Friedland~et.~al. \cite\refFGH\ have successfully generalized the classical Perron-Frobenius Theorem for nonnegative square matrices to nonnegative higher-order square tensors, based on different notions of irreducibility. However, via irreducibility in the sense of Definition \irreducible,\ the generalized Perron-Frobenius Theorem established in \cite \refCPZPF\ asserts the positive $H$-eigenvector is unique (up to a positive multiplier) in $\Bbb R_{+}^n\setminus\{0\}$, and there can be no $H$-eigenvectors on the boundary of the positive cone.  Whereas, via weak irreducibility in the sense of Definition \weakirrd,\ the generalized Perron-Frobenius Theorem in \cite\refFGH\ asserts that the positive $H$-eigenvector is unique (up to a positive multiplier) in $\Bbb R^n_{++}$, but it does not rule out the possibility of an $H$-eigenvector from the boundary of the positive cone.
More in-depth comparisons on the different kinds of spectral theory for tensors and their applications can be found in \cite\refCQZ.

When dealing with the $H$-eigenvalue problem for a nonnegative tensor $\Cal A$, the associated nonlinear map $\Cal T_{\Cal A}(x):=({\Cal A} x^{m-1})^{[{1\over {m-1}}]}: \Bbb R^n_+\rightarrow \Bbb R^n_+$ is indeed sub-linear, i.e. for any $0<t<1$, we have $t\Cal T_{\Cal A}(x)\le \Cal T_{\Cal A}(tx)$. Thus, using (weak) irreducibility of tensors, one may reduce the $H$-eigenvalue problem to a contraction mapping with respect to the Hilbert's projective metric defined in $\Bbb R^n_{++}$. It should be cautioned however, this particular technique is not applicable to the $Z$-eigenvalue problem, since the map $\Cal A x^{m-1}:  \Bbb R^n_+\rightarrow \Bbb R^n_+$ is not sub-linear for $m>2$. Aside from the purely algebraic differences between the two eigenvalue problems for tensors, the fundamentally different actions of the nonlinear maps induced by the nonnegative tensors in $\Bbb R^n_+$  described above also demand vastly different approaches; we refer to \cite\refCPZZ\ for a more systematic treatment of the general theory of $Z$-eigenvalue problem for nonnegative tensors.

\head{2. Basic Hypergraph Definitions}\endhead
\def\qctS{2}\qcts=0

We provide some definitions from the theory of hypergraphs.  The interested reader should refer to \cite{\refBerge}.
\setref\hypergraph
\definition{Definition \hypergraph}  Let $V$ be a finite set.  A hypergraph $H$ is a pair $(V,E)$, where $E \subseteq {\Cal P}(V)$, the power set of $V$.   The elements of $V$ are called the vertices and the elements of $E$ are called the edges.\enddefinition

We note that in the above definition of hypergraph, we do not allow for repeated vertices within an edge (often called a hyperloop).  In the case of hyperloops, we provide the following definition.  Recall the notion of a multiset is a generalization of the notion of a set in which members are allowed to appear more than once.

\setref\multihypergraph
\definition{Definition \multihypergraph}  Let $V$ be a finite set.  A multi-hypergraph $\Cal H$ is a pair $(V,E)$, where $E$ is a multi-set of subsets of $V$.   The elements of $V$ are called the vertices and the elements of $E$ are called the edges.\enddefinition

\setref\uniform
\definition{Definition \uniform}   A hypergraph $\Cal H$ is said to be $m$-uniform for an integer $m \ge 2$ if for all $e \in E$, the cardinal number of the subset $|e| =m$.  The term $m$-graph is often used in place of $m$-uniform hypergraph.  Similarly, a multi-hypergraph $\Cal H$ is said to be $m$-uniform for an integer $m \ge 2$ if for all $e \in E$, the cardinal number of the multi-set of $e$ is $m$.  Note that   the total number of elements in a multi-set, including repeated memberships, is the cardinality of the multi-set.  We will use the term $m$-multigraph in place of $m$-uniform multi-hypergraph.\enddefinition

\setref\completedef
\definition{Definition \completedef} Let $\Cal H=(V, E)$ be an $m$-graph. If any subset $A\subset V$ with $|A| =m$ is an edge of $\Cal H$, then $\Cal H$ is said to be $m$-complete.\enddefinition

\setref\regular
\definition{Definition \regular}  Let $\Cal H = (V, E)$ be a multi-hypergraph.  For each $v\in V$, the degree of $ v$ is the number of edges that contain it.  A hypergraph $\Cal H$ is said to be $r$-regular if every $v \in V$ has degree $r$.  \enddefinition

\setref\connected
\definition{Definition \connected}  In a multi-hypergraph $\Cal H=(V, E)$, a chain of length $q$ is a sequence $(x_1, E_1, x_2, E_2, \ldots, E_q, x_{q+1})$ such that 
\item\item{(1)}  $x_1, \ldots, x_q$ are all distinct vertices of $\Cal H$,
\item\item{(2)}  $ E_1, E_2, \ldots, E_q$ are all distinct edges of $\Cal H$, and
\item\item{(3)}  $x_k, x_{k+1} \in E_k$ for $k = 1, 2, \ldots, q$.

A multi-hypergraph $\Cal H$ is connected if there exists a chain starting at $x_1$ and terminating at $x_2$ for all $x_1, x_2 \in V$.  \enddefinition

\setref\partite
\definition{Definition \partite}  A hypergraph $\Cal H = (V,E)$ with $V=\{1, \ldots, n\}$ is said to be $k$-partite if there exists a partition of the vertices $V = V_1 \cup \cdots \cup V_k$ so that for any $k$ vertices $i_1, \ldots, i_k$ with $\{i_1, \ldots, i_k\} \in E$ each must lie in a distinct $V_i$ for $1 \le i \le k$. \enddefinition

\setref\adten
\definition{Definition \adten}  The adjacency tensor ${\Cal A}_{\Cal H}$ for an $m$-multigraph $\Cal H =(V,E)$ is the symmetric tensor ${\Cal A}_{\Cal H} = (a_{i_1 \ldots i_m}) \in {\Bbb R}^{[m,n]}$, where $n$ is the number of vertices and
$$a_{i_1 \ldots i_m} ={ 1\over (m-1)!} \cases 1, & \text{ if } i_1, \ldots, i_m \in E \cr
0, & \text{ otherwise.} \endcases$$
\enddefinition

\head {3. Main Results}\endhead
\def\qctS{3}\qcts=0

We begin this section by reviewing and proving some important results with regards to the $H$-eigenvalues of the adjacency tensor.  These results were previously obtained in \cite\refCooper\ regarding a general connected $m$-graph. We now present a short but different proof and show their results continue to hold for a general connected $m$-multigraph.  We begin with a lemma which verifies an $m$-mutligraph is connected if and only if the associated directed graph of its adjacency tensor is strongly connected.

\setref\lemmairred
\proclaim{Lemma \lemmairred}  Let $\Cal H$ be an $m$-multigraph on $n$ vertices.  Then $\Cal A_{\Cal H}\in \Bbb R^{[m,n]}_+$ is weakly irreducible if and only if $\Cal H$ is connected.  \endproclaim

\demo{Proof}  Assume $\Cal H$ is connected.   Recall the matrix $M(\Cal A_{\Cal H})=(m_{ij})$  where
$$ m_{ij}=\Sigma_{j\in \{i_2,\cdots, i_m \}} a_{ii_2,\cdots i_m } = |\{e\in E| i, j \in e\}|.$$
If $M(\Cal A_{\Cal H})$ is reducible, then there exists a nonempty proper $I \subset \{1, \ldots, n\}$ so that $m_{ij} =0$ for all $i \in I$ and all $j \not\in I$.  This implies $i, j$ are not in any edge for $i \in I$ and $j \not\in I$ and contradicts that $\Cal H$ is connected.  Using the equivalent notion following Definition \weakirrd,\ we see that since the matrix $M(\Cal A_{\Cal H})$ is irreducible; thus, $\Cal A_{\Cal H}\in \Bbb R^{[m,n]}_+$ is weakly irreducible. 

Assume $\Cal A_{\Cal H}\in \Bbb R^{[m,n]}_+$ is weakly irreducible.   Then its associated directed graph $G(\Cal A_{\Cal H})$ is strongly connected.   For every pair of vertices $x_1, x_{q+1}$, there is a directed path from $x_1$ to $x_{q+1}$.  This directed path consists of a sequence of edges $F_1, \ldots, F_q$ where $F_i$ is the directed edge from $x_i$ to $x_{i+1}$ in $G({\Cal H})$.  By definition, this implies there exists an edge $E_i$ in ${\Cal H}$ so that $E_i =\{x_i, x_{i+1},j_3, \ldots, j_m\}$ for some vertices $j_3, \ldots, j_m$.  We now have a sequence of edges $E_1, \ldots, E_q$ in ${\Cal H}$ and we verify $(x_1, E_1, x_2, E_2, \ldots, E_q, x_{q+1})$ reduces to a chain.  If $x_i$ are not distinct vertices of $\Cal H$, say $x_k = x_{k+m}$, then we may omit the piece of the chain from $E_k$ to $x_{k+m}$, inclusively.  If $E_i$ are not distinct edges, say $E_{k} = E_{k+m}$, then we may omit the piece of the chain from $x_{k+1}$ to $E_{k+m}$, inclusively.   After any necessary reduction, we have that $\Cal H$ is connected.  \qed\enddemo

\setref\Hpair
\proclaim{Theorem \Hpair} Let $\Cal H$ be a connected $m$-multigraph on $n$ vertices, then:
\item\item{1.}   There exists $(\lambda_0, x_0)\in {\Bbb R}_{++} \times {\Bbb R}_{++}^n$ an $H$-eigenpair of ${\Cal A}_{\Cal H}$, where 
$\lambda_0=\rho(\Cal A_{\Cal H})=\max_{x\in S_+} {\Cal A}_{\Cal H} x^m$ is the $H$-spectral radius of $\Cal A_{\Cal H}$, where $S_+=\{x\in \Bbb R^n_+\ |\ \sum_{i=1}^n x_i^m=1\}$.
\item\item{2.}  The positive eigenvector $x_0$ is unique (up to a positive multiplier) in $\Bbb R^n_{++}$.
\endproclaim
\demo{Proof} Since $\Cal H$ is connected, Lemma \lemmairred\ assures us that  $\Cal A_{\Cal H}\in \Bbb R^{[m,n]}_+$ is weakly irreducible. By Theorem 1.1 of \cite\refFGH,\ both the existence and uniqueness (up to a positive multiplier) of the positive $H$-eigenpair $(\lambda_0, x_0)\in {\Bbb R}_{++} \times {\Bbb R}_{++}^n$ follow readily.  The fact $\lambda_0=\rho(\Cal A_{\Cal H})$ follows from Corrollary 4.3 [12].  It remains to show $\lambda_0=\max_{x\in S_+} {\Cal A}_{\Cal H} x^m$.  Note the homogeneous polynomial ${\Cal A}_{\Cal H} x^m$ must achieve its absolute maximum value, denoted $\lambda_1$, on the compact set $S_+$, i.e. $\lambda_1=\max_{x\in S_+}{\Cal A}_{\Cal H} x^m$.  Since $\Cal A_{\Cal H}$ is a symmetric tensor, using Lagrange multipliers on the differentiable manifold $S_{++}$, the interior of $S_+$ (being an open subset of $\Bbb R^n$),  we see the eigenvector $x_0\in \Bbb R^n_{++}$ is a local maximizer of ${\Cal A}_{\Cal H} x^m$ on $S_+$; hence, we have  $\lambda_0\le \lambda_1$. Moreover, by the uniqueness of the $H$-eigenpair $(\lambda_0, x_0)\in\Bbb R_{++}\times  \Bbb R^n_{++}$, we see that 
$\lambda_0=\sup_{x\in S_{++}}{\Cal A}_{\Cal H} x^m$. Since the homogeneous polynomial ${\Cal A}_{\Cal H} x^m$ is nonnegative and continuous on $S_+$, we must have $$\sup_{x\in S_{++}}{\Cal A}_{\Cal H} x^m=\max_{x\in S_{+}}{\Cal A}_{\Cal H} x^m=\lambda_1,$$ which completes the proof.
\qed\enddemo

The direct approach for finding $\lambda_0$ via computing the multivariate resultant poses a computational challenge when the order $m$ and dimension $n$ both increase; in particular, if the order $m$ is high, then the problem is more difficult, if not impossible. However, it is worth pointing out, by making a small perturbation, the tensor $\Cal A_H+\mu \Cal E$ becomes positive, where $\mu>0$ is small and $\Cal E$ denotes the unit tensor whose entries consist of all ones.  According to \cite{\refCPZPrim, \refLZI, \refNQZ, \refQZ},\ starting with the initial vector $ \text{\bf 1}=(1, \cdots, 1)$, the NQZ-algorithm \cite\refNQZ,\ which is analogous to the power method, will converge linearly to the value $\rho(\Cal A_H+\mu\Cal E)$, where $\rho(\Cal A_H+\mu\Cal E)\rightarrow \rho(\Cal A_H)$ as $\mu\rightarrow 0$. For the sake of accessibility, we outline the NQZ-algorithm as follows:
\smallbreak
Given a nonnegative tensor $\Cal A\in \Bbb R_+^{[m,n]}$, not necessarily symmetric.
\smallbreak
Step 0.  Choose $x^{(0)} \in \Bbb R^n_{++}$.  Let $y^{(0)} = {\Cal A}(x^{(0)})^{m-1}$ and set $k :=0$.
\smallbreak
Step 1.  Compute $$\eqalign{
&x^{(k+1)} = {(y^{(k)})^{[{1\over m-1}]} \over \Vert (y^{(k)})^{[{1\over m-1}]} \Vert},\cr
&y^{(k+1)} = {\Cal A}(x^{(k+1)})^{m-1},\cr
&\underline\lambda_{k+1} = \min_{1\le i\le n}{(y^{(k+1)})_i \over (x_i^{(k+1)})^{m-1}},\cr
&\overline\lambda_{k+1} = \max_{1\le i\le n}{(y^{(k+1)})_i \over (x_i^{(k+1)})^{m-1}}.
}$$

\setref\exist
We now begin our studies on $Z$-eigenvalues of the adjacency tensor $\Cal A_{\Cal H}$ of a hypergraph $\Cal H$ by showing the existence of a positive $Z$-eigenvalue and a corresponding nonnegative $Z$-eigenvector. 
\proclaim{Theorem \exist}  Let $\Cal H$ be a an $m$-multigraph on $n$ vertices, then we have the following:
\item\item{1.}   There exists $(\lambda^*, x^*)\in {\Bbb R}_{++} \times {\Bbb R}_+^n$ a $Z$-eigenpair of ${\Cal A}_{\Cal H}$ such that
$$\lambda^*= \varrho(\Cal A_{\Cal H})=\max_{x \in S^{n-1}} {\Cal A}_{\Cal H} x^m=\max_{x \in S^{n-1}\cap {\Bbb R}_+^n} {\Cal A}_{\Cal H} x^m.$$
\item\item{2.} We have the following upper and lower bounds on $\lambda^*$:
$${ 1 \over  n^{m/2} (m-1)!}  \sum_{i\in V} \deg(i) \le \lambda^* \le \min\{D n^{1/2}, \ |E|\},$$ where $D$ is the maximal degree of any vertex in $V$ and $|E|$ denotes the cardinality of the edge set $E$. 
\endproclaim
\demo{Proof}  Since $\Cal A_{\Cal H}$ is a nonnegative symmetric $m$-order $n$-dimensional tensor which is not equal to the zero tensor, it is proven in \cite\refCPZZ  \ Corollary 3.12, that the largest real $Z$-eigenvalue $\lambda^*>0$ corresponds to the maximum value of the homogeneous polynomial  $\Cal A_{\Cal H} x^m$ on the standard unit sphere $S^{n-1}$ in $\Bbb R^n$. Furthermore, this value is attained at some point $x^*\in S^{n-1}\cap \Bbb R^n_+$.  Assertion 1 follows.

We consider the $m^{th}$ degree homogeneous polynomial ${\Cal A}_{\Cal H} x^m$  on $S^{n-1}$; its maximum value is greater than or equal to the value evaluated at $(1/\sqrt{n}, \ldots, 1/\sqrt{n})$, whose value is ${1 \over n^{m/2}  (m-1)!}  \sum_{i\in V} \text{deg}(i)$. This proves the first half of the chain of inequalities. As for the remaining inequality, since ${\Cal A}_{\Cal H} x^m$ has at most $|E|$ summands of the monomials $x_{i_1} x_{i_2}\cdots x_{i_m}$, each of which has maximum value $1$ on $S^{n-1}$, it is clear $\lambda^*\le |E|$. On the other hand, it follows from Proposition 3.3 of \cite\refCPZZ\ that 
$$\lambda^*\le \sqrt{n}\max_{1\le i\le n}{\sum_{i_2, \ldots, i_m=1}^n a_{i i_2\cdots i_m}}
=\sqrt{n}{\sum_{i_2, \ldots, i_m=1}^n a_{i_0 i_2\cdots i_m}}\le D\sqrt{n},$$
where $i_0 \in V$ is a vertex such that the degree of $i_0$ is $D$.  This completes our proof.\qed
\enddemo

We next investigate the condition needed to ensure the largest $Z$-eigenvalue $\lambda^*$ corresponds to a strictly positive $Z$-eigenvector. However, since the adjacency tensor $\Cal A_{\Cal H}$ is always reducible for $m$-graphs, such positivity result cannot be deduced directly from the existing literature. For this reason, we introduce another form of connectivity for $m$-graphs which is closely related to irreducible tensors, yet taking into account there are no hyperloops in an $m$-graph.  

\setref\redgraph
\definition{Definition \redgraph} An $m$-graph $\Cal H =(V, E)$ is not nicely-connected if there exists a nonempty $V_0 \subset V$ such that $|V_0|\le |V| -m +1$ and $\{i, i_2, \ldots, i_m\}\not\in E$ for all $i\in V_0$ and $i_2, \ldots, i_m \not\in V_0$.  Otherwise, we say $\Cal H$ is nicely-connected.  For an $m$-multigraph $\Cal H = (V, E)$, we say $\Cal H$ is not nicely-connected if there exists nonempty proper $V_0 \subset V$ such that $\{i, i_2, \ldots, i_m\}\not\in E$ for all $i\in V_0$ and $i_2, \ldots, i_m \not\in V_0$.    Otherwise, we say $\Cal H$ is nicely-connected.
\enddefinition

It follows directly from Definition \redgraph\ that in order for either an $m$-multigraph or an $m$-graph to be nicely-connected, it must be connected a priori.  When $m=2$, a graph is nicely-connected is equivalent to its being connected and the irreducibility of its adjacency matrix.  For $m >2$, a multi-graph is nicely-connected is equivalent to the irreducibility of the adjacency tensor.  However, it is possible to have a nicely-connected $m$-graph whose adjacency tensor is reducible in the sense of Definition \irreducible\ and weakly irreducible in the sense of Definition \weakirrd,\ see Example \exampleD.    Furthermore, the nicely-connectedness is intrinsically determined by a hypergraph, not by any tensor associated with it.

\setref\compa
\example{Example \compa} Consider the $3$-graph $\Cal H=(V,E)$, where $V=\{1,2,3,4,5\}$ and $E=\{123, 345\}$.   Then $\Cal H$ is not a regular $3$-graph, but $\Cal H$ is nicely-connected. This example will be reexamined as Example \exampleD.
\endexample

\setref\compb
\example{Example \compb} Consider the $3$-graph $\Cal H=(V,E)$, where $V=\{1,2,3, 4,5, 6\}$ and $E=\{123, 236, 145, 456\}$.   Then $\Cal H$ is a $2$-regular connected $3$-graph, but $\Cal H$ is not nicely-connected as seen by taking $V_0 = \{4, 5\}$.
\endexample

In \S 4, we will also demonstrate via a specific example (see Example \exampleEe)\ that a connected regular uniform multigraph is not necessarily nicely-connected. Thus, our notion of nicely-connected is indeed stronger than connectivity but structurally different from regularity.

\setref\positive
\proclaim{Theorem \positive} Let $\Cal H$ be an $m$-multigraph on $n$ vertices.  If $\Cal H$ is nicely-connected,  then any $Z$-eigenpair of $\Cal A_{\Cal H}$ such that $(\lambda, x)\in {\Bbb R}_{++} \times {\Bbb R}_+^n$   must satisfy $x_i >0$ for all $1\le i\le n$.  Conversely, if $\Cal H$ is not nicely-connected, then there exists a positive $Z$-eigenvalue which corresponds to an eigenvector whose coordinates are all nonnegative but at least one coordinate is zero.
\endproclaim

\demo{Proof}  By the first assertion of Theorem \exist,\ the set of $Z$-eigenpairs $(\lambda, x)$ of the adjacency tensor $\Cal A_{\Cal H}$ such that $(\lambda, x)\in {\Bbb R}_{++} \times {\Bbb R}_+^n$ is nonempty. Let $(\lambda, x)\in {\Bbb R}_{++} \times {\Bbb R}_+^n$ be a $Z$-eigenpair of $\Cal A_{\Cal H}$. Let $V_0$ be the set of indices for which $x_i =0$.  

Let  $E':=\{e| e\in E, e\subset V\setminus V_0\}$.  Let $V'=\{v\in V| v\in e \text{ for some } e \in E'\}$.  We consider the $m$-multigraph $\Cal H'=(V', E')$.  

Suppose $V'=\emptyset$.   Then all edges of $E$ contain a vertex from $V_0$.  If $i \not\in V_0$, then the equation $({\Cal A}_{\Cal H} x^{m-1})_i \ne 0$ is a contradiction.  Hence, $V' \ne \emptyset$.   For a multigraph, it is sufficient for $V'$ to be nonempty.  For an $m$-graph, we must also make certain the set $V_0$ is not too large. Suppose $|V_0| >n -m+1$ and $\Cal H$ is an $m$-graph (i.e. there are no hyperloops).  Then all edges of $E$ contain a vertex from $V_0$.  If $i \not\in V_0$, then the equation $({\Cal A}_{\Cal H} x^{m-1})_i \ne 0$ is a contradiction.  So, for any $m$-graph, we have $|V_0| \le n -m+1$.

Let $j \in V\setminus V'$.  Then the only edges that contain $j$  also contain a vertex from $V_0$, which implies the equation $(\Cal A_{\Cal H} x^{m-1})_j =0 $; hence, $j \in V_0$.  Therefore, $V= V_0\ \dot\cup\ V'$.  This shows $\Cal H$ is not nicely-connected since $\{i, i_2, \ldots,i_m\} \not\in E$ for $i \in V_0$ and $\{i_2, \ldots, i_m\} \subset V'$.  Contrapositively, the conclusion follows.

For the second statement, we assume $\Cal H$ is not nicely-connected with nonempty proper $V_0 \subset V$ so that $\{i, i_2, \ldots, i_m\} \not\in E$ for all $i \in V_0$ and all $i_2, \ldots, i_m \not\in V_0$ with $|V_0| \le n-m+1$ if $\Cal H$ is an $m$-graph.  Let $V'' = V\setminus V_0$ and $E'':=\{e\in E| e\subset \Cal P(V'')\}$.   Consider the hypergraph ${\Cal H}'' =(V'', E'')$.   By Theorem \exist, there exists $(\lambda^*, x^*)\in {\Bbb R}_{++} \times {\Bbb R}_+^{|V\setminus V_0|}$ a $Z$-eigenpair of ${\Cal A}_{\Cal H''}$.  Let $y^*\in {\Bbb R}^n$ be defined via $y_i = (x^*)_i$ if $i \in V\setminus V_0$ and $y_i = 0$ if $i \in V_0$.  Then we claim that $(\lambda^*, y)$ is the desired $Z$-eigenpair of $\Cal H$.   For $i \in V$, we have $({\Cal A}_{\Cal H''} {x^*}^{m-1})_i =\lambda^* (x^*)_i = \lambda^* y_i$.  We may add to the left hand side of this equation $0 = \sum_{}a_{ii_2, \ldots i_m} y_{i_2} \cdots y_{i_m}$, where the sum is taken over all $i_2, \ldots i_m$ for which $\{i, i_2, \ldots, i_m\} \in E$ and at least one of the indices from $\{i_2, \ldots, i_m\}$ is in $V_0$ (so at least one factor of each monomial is zero).  We then have $({\Cal A}_{\Cal H} {y}^{m-1})_i =\lambda^* y_i$.  Suppose $i \in V_0$.  Then $y_i =0$, so it suffices to show $({\Cal A}_{\Cal H} {y}^{m-1})_i= 0$.  Since $i \in V_0$, any edge that contains $i$ must contain another vertex from $V_0$; hence, $({\Cal A}_{\Cal H} {y}^{m-1})_i= 0$ as required to verify $(\lambda^*, y)$ is the desired $Z$-eigenpair of $\Cal H$.  \qed\enddemo
\smallbreak\noindent
{\bf Remark.} We comment that based on the above Theorem \positive,\ the existence of a strictly positive $Z$-eigenvector $x^*$ enables us to detect whether an $m$-multigraph is nicely-connected. This feature depicts a more subtle combinatorial structure of a connected $m$-graph than the $H$-eigenvalues encodes.

\setref\complete
\proclaim{Corollary \complete} If $\Cal H$ is an $m$-graph on $n$ vertices which is $m$-complete,  then any $Z$-eigenpair $(\lambda, x)\in {\Bbb R}_{++} \times {\Bbb R}_+^n$ of $\Cal A_{\Cal H}$  must satisfy $x_i >0$ for all $1\le i\le n$.  
\endproclaim

\demo{Proof}   It is sufficient to show that a complete $m$-graph  $\Cal H$ is nicely connected.  Suppose not.  Then there must be a nonempty proper subset $V_0\subset V$ such that $|V_0| \le |V| -m+1$ and $\{i, i_2, \ldots, i_m\} \not\in E$ for all $i\in V_0$ and $i_2, \ldots i_m \not\in V_0$. But $\Cal H$ is complete, so $\{i, i_2, \cdots, i_m\}\in E$; this is a contradiction.  \qed\enddemo
\setref\compute
\proclaim{Corollary \compute}  If $\Cal H$ is a connected $m$-multigraph on $n$ vertices, then the algorithm (SS-HOPM cf. \cite\refKM) can be used to compute a $Z$-eigenvalue of $\Cal A_{\Cal H}$ as follows:\endproclaim
\smallbreak
Given a symmetric tensor $\Cal A\in \Bbb R^{[m,n]}$.
\smallbreak
Step 0. Choose  $x_{(0)}\in \Bbb R^n_+\setminus \{0\}$, set $\lambda_0=\Cal A x_{(0)}^m$, and
choose a sufficiently large shift constant $\alpha$, e.g. $$\alpha=\lceil m\sum_{i_1,\cdots, i_m=1}^n a_{i_1\ldots i_m}\rceil,$$
where $\lceil \gamma\rceil$ is the ceiling function, i.e. it equals the smallest integer  no less than $\gamma$.
\smallbreak
Set $k:=0$.
\smallbreak
Step 1. Set $y_{(k+1)}:=\Cal A x_{(k)}^{m-1}+\alpha x_{(k)}$.
\smallbreak
Step 2. Compute 
$x_{(k+1)}:={y_{(k+1)}\over \Vert y_{(k+1)}\Vert}$ and 
$\lambda_{k+1}:=\Cal A x_{(k+1)}^m$.

\setref\largestgraph
\proclaim{Corollary \largestgraph}  If $\Cal H$ is an $r$-regular $m$-graph on $n$ vertices, then $\lambda_0=r n^{-(m-2)/2}$ is a $Z$-eigenvalue with corresponding positive eigenvector $x_0:={{\text{\bf 1}}/\sqrt{n}}$, where ${\text{\bf 1}}= (1, \ldots, 1)$.   \endproclaim
 \demo{Proof}  We notice that $r$-regularity of an $m$-graph implies that $(\Cal A_{\Cal H}x^{m-1})_i$ has $r$ summands of monomials with degree $m-1$.  We now verify that $x_0$ is an eigenvector of $\lambda_0$:
  $$\bigg({\Cal A_{\Cal H}} { x_0}^{m-1}\bigg)_i = r \cdot (1/\sqrt{n})^{m-1}= (r(1/\sqrt{n})^{m-2}) (x_0)_i\quad  \text{for\ all}\ 1\le i\le n.\ \ \qed$$\enddemo

\setref\largestcomplete
\proclaim{Corollary \largestcomplete} If $\Cal H$ is an $m$-graph on $n$ vertices which is $m$-complete, then $\lambda_0= {n-1 \choose m-1} n^{-(m-2)/2}$ is a $Z$-eigenvalue with a corresponding positive eigenvector. \endproclaim

\demo{Proof}  Since $\Cal H$ is an $m$-graph on $n$ vertices which is $m$-complete, we see that each vertex is contained in ${n-1 \choose m-1}$ edges; hence, $\Cal H$ is ${n-1 \choose m-1}$-regular.  \qed\enddemo

\setref\symmetric
\proclaim{Theorem \symmetric}  If $\Cal H$ is an $m$-graph on $n$ vertices, then 
\item\item{1.}  If $m$ is odd and $\lambda $ is an $E$-eigenvalue of $\Cal A_{\Cal H}$, then $-\lambda$ is an $E$-eigenvalue of $\Cal A_{\Cal H}$.
\item\item{2.}  If $\Cal H$ is $m$-partite with $m$ even and $\lambda $ is an $E$-eigenvalue of $\Cal A_{\Cal H}$, then $-\lambda$ is an $E$-eigenvalue of $\Cal A_{\Cal H}$.

In either case above, the sum of all $E$-eigenvalues is zero.
\endproclaim
\demo{Proof}  If $m$ is odd and $(\lambda, x)$ is a $E$-eigenpair, then $(-\lambda, - x)$ is also a $E$-eigenpair since ${\Cal A_{\Cal H}}(-x)^{m-1} = {\Cal A_{\Cal H}}(x)^{m-1}$ since $m$ is odd and this implies $${\Cal A_{\Cal H}}(-x)^{m-1} = {\Cal A_{\Cal H}}(x)^{m-1} = \lambda x = (-\lambda)(-x).$$ 
 
 Suppose $\Cal H$ is an $m$-partite $m$-graph with $m$ even.  Then there exists a partition $V = V_1 \cup \cdots \cup V_m$ so that $a_{i_1 \cdots i_m} \ne 0$ implies $i_j$ must each lie in a distinct $V_{i_j}$.  Suppose $\lambda$ is a $E$-eigenvalue with corresponding eigenvector $ x$.  Let $ y$ be defined via $$y_i = \cases -x_i & \text{if } i \in V_1 \cr x_i& \text {if } i \not\in V_1. \endcases$$  We then consider $({\Cal A_{\Cal H}}y^{m-1})_i$.  If $i \in V_1$, then $({\Cal A_{\Cal H}}y^{m-1})_i$ has no terms that contain a $y_1$ factor; hence, $$({\Cal A_{\Cal H}}y^{m-1})_i=({\Cal A_{\Cal H}}x^{m-1})_i = \lambda x_i = (-\lambda)(-x_i) = (-\lambda) y_i.$$  If $i \not\in V_1$, then every term of $({\Cal A_{\Cal H}}y^{m-1})_i$ has exactly one factor $y_{i_j} \in V_1$.  Hence, $$({\Cal A_{\Cal H}}y^{m-1})_i= -({\Cal A_{\Cal H}}x^{m-1})_i=-\lambda x_i = -\lambda y_i.\qed$$   \enddemo

\head {4. Examples and Open Problems}\endhead
\def\qctS{4}\qcts=0

\setref\exampleAa
\example{Example \exampleAa}  We now give an example of a connected 3-graph that is not nicely-connected to illustrate the failure of the existence of a positive $Z$-eigenvector, demonstrating the intricacies of Theorem \positive.\ Let $\Cal H$ be the hypergraph given by $V=\{1,2,3,4,5,6,7\}$ and $E=\{123, 345, 567\}$. The corresponding $Z$-eigenvalue problem of the adjacency tensor is given by 
$$\eqalign{&x_2 x_3= \lambda x_1\cr
&x_1 x_3 = \lambda x_2\cr
&x_1 x_2+x_4 x_5= \lambda x_3\cr
&x_3 x_5 = \lambda x_4\cr
&x_3 x_4+x_6 x_7= \lambda x_5\cr
&x_5 x_7= \lambda x_6\cr
&x_5 x_6 = \lambda x_7\cr
&x_1^2+x_2^2+x_3^2+x_4^2+x_5^2+x_6^2+x_7^2 = 1}.$$
By computing the elimination ideal, we see that the $E$-characteristic polynomial of $\Cal A_{\Cal H}$ is  $-3\lambda^3+\lambda$, whose three roots are $0, \pm {{\sqrt 3}\over 3}$. However, the largest $Z$-eigenvalue $\lambda^*={{\sqrt 3}\over 3}$ does not correspond to any strictly positive eigenvectors.  The $Z$-eigenvalue $\lambda^*={{\sqrt 3}\over 3}$ has infinitely many nonnegative $Z$-eigenvectors given by 
\par\noindent $(0,0, t, t, {\sqrt{3}\over 3},  {\sqrt{1-3t^2} \over \sqrt{3}}, {\sqrt{1-3t^2} \over \sqrt{3}})$ and $({\sqrt{1-3t^2} \over \sqrt{3}}, {\sqrt{1-3t^2} \over \sqrt{3}}, {\sqrt{3} \over 3}, t, t, 0, 0)$ for $-{\sqrt{3} \over 3} \le t\le {\sqrt{3} \over 3}$.  Then by considering positive $Z$-eigenvalues with nonnegative $Z$-eigenvectors, we have $V_0=\{ 1,2\}$ if $t \ne  \pm {\sqrt{3} \over 3}$ or  $V_0=\{1,2,6,7\}$ if $t = \pm {\sqrt{3} \over 3}$,  showing that $\Cal H$ is not nicely-connected.  
\endexample

\setref\exampleC
\example{Example \exampleC}  We give an example of a connected, yet not nicely-connected, $3$-multigraph where the eigenvector corresponding to a positive real eigenvalue is not positive in all coordinates.  Let $V=\{1,2\}$ and $E=\{112, 222\}$.  Then the  corresponding $Z$-eigenvalue problem of the adjacency tensor is given by 
$$\eqalign{&x_1x_2=\lambda x_1\cr
&{1\over 2} x_1^2 + {1 \over 2} x_2^2 = \lambda x_2\cr
&x_1^2+x_2^2=1}$$
The $E$-eigenvalues are $\pm {1\over 2}, \pm {\sqrt{2}\over 2}$.  The eigenvector corresponding to ${1 \over 2}$ is $(0,1)$.  However, the eigenvector corresponding to ${\sqrt{2}\over 2}$ is $({\sqrt{2} \over 2}, {\sqrt{2}\over 2})$.  
\endexample

\setref\exampleEe
\example{Example \exampleEe}  We give an example of a connected, yet not nicely-connected, $2$-regular $3$-multigraph.  Let $V=\{1,2\}$ and $E=\{111, 112, 222\}$.  Then the corresponding $Z$-eigenvalue problem of the adjacency tensor is given by
$$\eqalign{& {1\over 2} x_1^2 + x_1x_2 = \lambda x_1\cr
& {1 \over 2} x_1^2 + {1\over 2} x_2^2 = \lambda x_2\cr
&x_1^2+x_2^2 = 1}$$
The largest $Z$-eigenvalue of $0.951057$ corresponds to a strictly positive $Z$-eigenvector of $(0.850651,0.525731)$. However,  the positive $Z$-eigenvalue of $\lambda=0.5$ corresponds to a $Z$-eigenvector of $(0,1)$; it is this second $Z$-eigenpair that demonstrates the failure of the graph to be nicely-connected with $V_0 = \{1\}$.  \endexample
\setref\exampleD
\example{Example \exampleD}  We give an example of a nicely-connected non-regular $3$-graph where there are infinitely many positive $Z$-eigenvectors which correspond to the largest positive $Z$-eigenvalue.   Let $V=\{1,2,3, 4, 5\}$ and $E=\{123,345\}$.  Then the  corresponding $Z$-eigenvalue problem of the adjacency tensor is given by 
$$\eqalign{&x_2x_3=\lambda x_1\cr
& x_1x_3= \lambda x_2\cr
&x_1x_2+x_4x_5= \lambda x_3\cr
&x_3x_5= \lambda x_4\cr
&x_3x_4= \lambda x_5\cr
&x_1^2+x_2^2+x_3^2+x_4^2+x_5^2 =1}.$$  The $Z$-eigenvalues are $0, \pm {\sqrt{3}\over 3}$.  For the largest $Z$-eigenvalue $\lambda^* = {\sqrt{3} \over 3}$, we calculate its corresponding positive $Z$-eigenvectors to be given by $x_3 ={\sqrt{3}\over 3}$, $x_1^2+x_5^2 = {1 \over 3}$, $x_2=x_1>0$, and $x_5=x_4>0$.   Hence, we note that the largest $Z$-eigenvalue of a connected $m$-hypergraph may have more than one $Z$-eigenvector (i.e. it is not real geometrically simple).\endexample

\setref\exampleB
\example{Example \exampleB}  We give an example of a 1-regular 4-graph which is 4-partite that has symmetric $E$ spectrum, illustrating Theorem \symmetric.  Let $\Cal H$ be the hypergraph given by $V=\{1,2,3,4\}$ and $E=\{1234\}$.  Then the corresponding $E$-eigenvalues of the adjacency tensor are $0, \pm {1\over 4}$.  

The $E$-eigenvectors of $0$ are:  
\par\qquad $( 0, 0, t, -\sqrt{1 - t^2});\ (0, 
   0, t, \sqrt{1 - t^2});$
   \par\qquad $( 0, t, 0, -\sqrt{1 - t^2});\ ( 0, t, 0,\sqrt{1 - t^2}); $
  \par\qquad $(0, t,
  -\sqrt{1 - t^2}, 0);\ (0, t,
  \sqrt{1 -t^2}, 0);$
  \par\qquad $ (t, 0, 0, 
   -\sqrt{1 - t^2});\ (t, 0, 0, 
   \sqrt{1 - t^2});$
   \par\qquad $ (t, 0, 
   -\sqrt{1 - t^2}, 0);\ (t,  0, 
   \sqrt{1 -t^2}, 0);
   $\par\qquad $(t, -\sqrt{1 - t^2},0,0);\ (t,    \sqrt{1 - t^2},0,0)$.
   
   The $E$-eigenvectors of ${1 \over 4}$ are the eight vectors coming from  $(\pm {1 \over 2}, \pm {1 \over 2}, \pm {1 \over 2}, \pm {1 \over 2})$, where either all are negative, all are positive, or there are exactly two positive and two negatives.
   
   The $E$-eigenvectors of $-{1 \over 4}$ are the eight vectors coming from $(\pm {1 \over 2}, \pm {1 \over 2}, \pm {1 \over 2}, \pm {1 \over 2})$, where either there are exactly three negatives or exactly one negative.
  \endexample

\setref\exampleA
\example{Example \exampleA}  We give an example of a 3-regular 3-graph that has  symmetric $Z$-spectrum, illustrating Corollary \largestgraph.  This is a complete $3$-graph on $4$ vertices and verifies the conclusion of Corollary 3.11.  We note that the positive $Z$-eigenvalues that do not have a positive $Z$-eigenvector do not violate Theorem 3.7.  Let $\Cal H$ be the hypergraph given by $V=\{1,2,3,4\}$ and $E=\{ 123, 124, 134, 234\}$.  Then the corresponding $Z$-eigenvalue problem of the adjacency tensor is given by 
$$\eqalign{&x_2 x_3+x_2x_4+x_3x_4 = \lambda x_1\cr
&x_1 x_3+x_1x_4+x_3x_4 = \lambda x_2\cr
&x_1 x_2+x_1x_4+x_2x_4 = \lambda x_3\cr
&x_1 x_2+x_1x_3+x_2x_3 = \lambda x_4\cr
&x_1^2+x_2^2+x_3^2+x_4^2 = 1}.$$
The $Z$-eigenvectors of each $Z$-eigenvalue are:

For $\lambda=0$, we have $(\pm1, 0, 0, 0),(0, \pm 1, 0, 0 ,0),(0,0,\pm1,0),(0,0,0, \pm1)$. 

For $\lambda = {3\over 2}$, we have $ ({1\over 2}, {1\over 2},{1\over 2},{1\over 2})$. 

For $\lambda =-{3 \over 2}$, we have $(-{1\over 2}, -{1\over 2},-{1\over 2},-{1\over 2})$.  

For $\lambda = {\sqrt{6} \over 6}$, we have the following $Z$-eigenvectors:
\par\qquad $(-0.66056, 0.252311, -0.66056,  0.252311);$
\par\qquad  $(0.252311,  -0.66056,0.252311, -0.66056);$
\par\qquad  $(-0.66056,  0.252311,  0.252311-0.66056);$
\par\qquad  $(0.252311, 0.252311, -0.66056,  -0.66056);$
\par\qquad  $(-0.66056,  -0.66056,  0.252311, 0.252311);$
\par\qquad   $(0.252311,  -0.66056, -0.66056, 0.252311)$. 

For $\lambda = -{\sqrt{6} \over 6}$, we have the opposite of the eigenvectors given above for ${\sqrt{6} \over 6}$.  

For $\lambda = { 4 \sqrt{21} \over 21}$, we have the following $Z$-eigenvectors:
\par\qquad $( -0.436436,  0.654654, 
 -0.436436,  -0.436436);$
 \par\qquad$ ( -0.436436,  -0.436436,  -0.436436, 
  0.654654);$
  \par\qquad $( 0.654654, 
  -0.436436,  -0.436436, 
  -0.436436);$
  \par\qquad $(-0.436436, 
 -0.436436,  0.654654, -0.436436)$.

For $\lambda=-{ 4 \sqrt{21} \over 21}$, we have opposite of the eigenvectors given above for ${ 4 \sqrt{21} \over 21}$.
\endexample

\setref\examplecomplete
\example{Example \examplecomplete}  Let $\Cal H$ be the complete $3$-graph on $5$ vertices.  We then compute via \cocoa \cite{\refCocoa}\ that an eigenvalue is a root of the polynomial $-6360\lambda^9 + 56507\lambda^7 - 81513\lambda^5 + 31833\lambda^3 - 2916\lambda$.  Hence, the eigenvalues are $0, \pm \sqrt{{1 \over 318}(89-13\sqrt{13})}, \pm \sqrt{{89 \over 318} + {13\sqrt{13} \over 318}}, \pm {3 \over 2\sqrt{2}}, \pm {6 \over \sqrt{5}}$.  We note the largest positive eigenvalue has only one eigenvector.   The other positive eigenvalues each has more than one eigenvector. \endexample

\setref\questionA
\example{Question \questionA}  In an $m$-multigraph $\Cal H$, what conditions are needed such that $\Cal A_{\Cal H}$ has a unique positive $Z$-eigenvalue?  \endexample

\setref\questionB
\example{Question \questionB}  In a connected regular $m$-graph, does the largest $Z$-eigenvalue have only one corresponding $Z$-eigenvector (i.e. is it real geometrically simple)?\endexample

\setref\questionC
\example{Question \questionC}  In a connected $m$-graph, what conditions are needed such that $\Cal A_{\Cal H}$ has symmetric spectrum in the sense of Definition \Hpairdef?
\endexample


\Refs \widestnumber\key{99}

\ref\refkey{\refBerge} \by Claude Berge \book Graphs and hypergraphs \bookinfo North-Holland Mathematical Library, second edition \vol 6 \publ North-Holland\yr 1976 \endref

\ref\refkey{\refBP} \by S.R.~Bul\`{o} and M. Pelillo \paper  New bounds on the clique
number of graphs based on spectral hypergraph theory, in: T.
St\"utzle ed., Learning and Intelligent Optimization, Springer Verlag, Berlin, \yr 2009 \pages 45-48 \endref

\ref\refkey{\refCS} \by D.~Cartwright and B.~Sturmfels \paper The number of eigenvalues of a tensor  \jour Linear Algebra Appl. (in press)\endref

\ref\refkey{\refCPZ} \by K.C.~Chang, K.~Pearson, and T.~Zhang \paper On eigenvalue problems of real symmetric tensors  \jour  J. Math. Anal. Appl. \yr 2009 \vol 350 \pages 416-422\endref

\ref\refkey{\refCPZPF} \by K.C.~Chang, K.~Pearson, and T.~Zhang \paper Perron-Frobenius theorem for nonnegative tensors \jour Commun. Math. Sci. \yr 2008\vol 6  \issue 2 \pages 507-520\endref

\ref\refkey{\refCPZPrim} \by K.C.~Chang, K. ~Pearson, and T.~Zhang \paper Primitivity, the convergence of the NQZ method, and the largest eigenvalue for nonnegative tensors \jour SIAM. J. Matrix Anal. \& Appl., \yr 2011 \vol 32 \pages 806-819\endref

\ref\refkey{\refCPZZ} \by K.C.~Chang, K.~Pearson and T.~Zhang \paper Some variational principles of the $Z$-eigenvalues for nonnegative tensors \jour School of
Mathematical Sciences, Peking University, December \yr 2011\endref

\ref\refkey{\refCQZ} \by K.C.~Chang, L.~Qi, and T.~Zhang \paper A survey on the spectral theory of nonnegative tensors, preprint\endref

\ref\refkey{\refChung} \by F.R.K.~Chung\book Spectral graph theory\yr 1997\bookinfo Am. Math. Soc. \publ Providence\endref

\ref\refkey{\refCocoa} \by {{CoCoA}Team} \book 
    {{{\hbox{\rm C\kern-.13em o\kern-.07em
C\kern-.13em o\kern-.15em A}}}: a system for doing {C}omputations
in {C}ommutative {A}lgebra} \publ
   {Available at \/ {\tt http://cocoa.dima.unige.it}}\endref
   

\ref\refkey{\refCooper} \by J.~Cooper and A.~Dutle \paper Spectra of Uniform Hypergraphs,
Department of Mathematics, University of South Carolina, June 2011,
arXiv:1106.4856, 2011\endref



\ref\refkey{\refFGH}\by S.~Friedland, S.~Gaubert, and L.~Han \paper Perron-Frobenius theorem for nonnegative multilinear forms and extensions \jour Linear Algebra Appl. (in press) \endref

\ref\refkey{\refHQ} \by S.~Hu and L.~Qi \paper Algebraic connectivity of an even uniform
hypergraph \jour J. Comb. Optim. (in press) \endref

\ref\refkey{\refHu} \by S.~Hu \paper A note on the positive eigenvector of nonnegative tensors\yr 2010 preprint\endref 

\ref\refkey{\refKM} \by T.~Kolda and J.~Mayo \paper Shifted Power method for computing tensor eigenpairs \jour SIAM J. Matrix Anal. Appl., \yr 2011 \vol 34
\pages1095-1124\endref

\ref\refkey{\refLima} \by L.H.~Lim \paper Singular values and eigenvalues of tensors, A
variational approach\jour Proc. 1st IEEE International workshop on
computational advances of multi-tensor adaptive processing, Dec.
13-15\yr 2005\pages 129--132\endref

\ref\refkey{\refLimb} \by L.H.~Lim \paper Multilinear pagerank: measuring
higher order connectivity in linked objects \jour The Internet : Today
and Tomorrow \vol July \yr 2005 \endref

\ref\refkey{\refLZI} \by Y.~Liu, G.~Zhou, and N.F.~Ibrahim \paper An always convergent algorithm for the largest eigenvalue of an irreducible nonnegative tensor  \jour J. Comput. Appl. Math.  \yr 2010\vol 235 \pages 286-292\endref

\ref\refkey{\refNQZ}\by M.~Ng, L.~Qi, and G.~Zhou \paper Finding the largest eigenvalue of a nonnegative tensor \jour SIAM. J. Matrix Anal. \& Appl. \yr 2009\vol  31, \issue 3 \pages 1090-1099\endref



\ref\refkey{\refQi} \by L.~Qi  \paper Eigenvalues of a real supersymmetric tensor\jour J. Symbolic Comput. \yr 2005\vol 40 \pages 1302--1324 \endref

\ref\refkey{\refQZ} \by L.~Qi and L.~Zhang \paper Linear convergence of an algorithm for computing the largest eigenvalue of a nonnegative tensor \jour Numer. Linear Algebra Appl. (in press)\endref


\endRefs
\enddocument